\documentclass[12pt,oneside]{article}

\usepackage{amsfonts}
\usepackage{amscd}
\usepackage{amsthm}
\usepackage{amssymb}
\usepackage{amsmath}
\usepackage{epsfig}

\usepackage{graphicx}
\usepackage{xypic}
\input xypic
\input xy
\xyoption{all}

\title{EZ-structures and topological applications.}

\author{F. T. Farrell \footnote{This research was supported in part
by the National Science Foundation.} \hskip 5pt \& J.-F. Lafont}

\theoremstyle{definition}
\newtheorem{Def}{Definition}[section]

\theoremstyle{proposition}
\newtheorem{Lem}{Lemma}[section]
\newtheorem{Prop}{Proposition}[section]
\newtheorem{Claim}{Claim}
\newtheorem{Claim2}{Claim}
\newtheorem{Fact}{Fact}

\theoremstyle{plain}

\newtheorem{Thm}{Theorem}[section]

\theoremstyle{remark}

\newtheorem*{Prf}{Proof}

\begin{document}

\maketitle

\begin{abstract}
In this paper, we introduce the notion of an EZ-structure on a group, an equivariant 
version of the Z-structures introduced by Bestvina [4].  Examples of groups having an
EZ-structure include:
\begin{itemize}
\item torsion free $\delta$-hyperbolic groups.
\item torsion free $CAT(0)$-groups.
\end{itemize}
Condition (*) was introduced by Farrell-Hsiang [8] in order to provide an 
abstract setting in which to prove the Novikov conjecture.  We introduce a generalization,
termed condition $(*_\Delta)$ of condition $(*)$, and show that
any group that has an EZ-structure automatically satisfies condition $(*_\Delta)$.  The 
argument of Farrell-Hsiang extends to show that the Novikov conjecture holds for any
group satisfying condition $(*_\Delta)$.  As another application of these techniques,
we show how, in the case of a $\delta$-hyperbolic group $\Gamma$, we can obtain a lower 
bound for the homotopy groups $\pi_n(\mathcal P(B\Gamma))$, where $\mathcal P(\cdot )$ is
the stable topological pseudo-isotopy functor.
\end{abstract}

\section{Introduction.}

Let $\Gamma$ be a discrete group.  Bestvina [4] defined the notion of a Z-structure
on $\Gamma$ as a pair $(\bar X,Z)$ of spaces satisfying the following four axioms:

\begin{itemize}
\item $\bar X$ is a Euclidean retract (ER).
\item $Z$ is a Z-set in $\bar X$.
\item $\bar X-Z$ admits a fixed point free, properly discontinuous, cocompact 
action by the group $\Gamma$.
\item The collection of translates of a compact set in $\bar X-Z$ forms a null sequence in 
$\bar X$; i.e. for every open cover $\mathcal{U}$ of $\bar X$, all but finitely
many translates are $\mathcal{U}$ small.
\end{itemize}

Let us now introduce an equivariant version of a Z-structure:

\begin{Def}
We say that $(\bar X,Z)$ is an EZ-structure (equivariant Z-structure) on $\Gamma$ 
provided that $(\bar X,Z)$ is a Z-structure, and in addition, the $\Gamma$ action on
$\bar X-Z$ extends to an action on $\bar X$.
\end{Def}

Examples of groups with an EZ-structure include torsion-free $\delta$-hyperbolic groups [3]
and $CAT(0)$-groups [4].  We note that a special case of a Z-structure on $\Gamma$
is the situation where $\bar X$ is a disk $\mathbb{D}^n$, and $Z=\partial \mathbb{D}^n
=S^{n-1}$:

\begin{Def}
We say that $\Gamma$ satisfies condition (*) provided that there is an EZ-structure
of the form $(\mathbb{D}^n, S^{n-1})$.
\end{Def}

Farrell-Hsiang introduced this special case in [8] (see also [9],
[12], [13]).  Their motivation for the development of condition (*) was
that it provided an abstract setting under which the Novikov conjecture could be 
verified for the group $\Gamma$.  Observe that there are groups with an EZ-structure 
that do 
{\it not} satisfy condition (*);  for example, the free group on 2-generators.  
We now introduce a condition $(*_\Delta)$ for torsion-free groups, 
generalizing condition (*).  (For non torsion-free groups see Definition 3.1 below)

\begin{Def}
We say that $\Gamma$ satisfies condition $(*_\Delta)$ provided that there is an EZ-structure
of the form $(\mathbb{D}^n, \Delta)$, where $\Delta$ is a closed subset of $\partial
\mathbb{D}^n=S^{n-1}$
\end{Def}

We are now ready to state the first two theorems of this paper:

\begin{Thm}
Let $\Gamma$ be a discrete group, and assume that $\Gamma$ has an EZ-structure.  Then
$\Gamma$ satisfies condition $(*_\Delta)$.
\end{Thm}

\begin{Thm}
Let $\Gamma$ be a torsion-free discrete group satisfying condition $(*_\Delta)$.  
Then the Novikov conjecture holds for the group $\Gamma$.
\end{Thm}

The proofs of these theorems will be provided in section 2 and section 3 respectively.  
We note that the second theorem
is not new, as Carlsson-Pederson [6] have already proven that groups with an EZ-structure
satisfy this form of the Novikov conjecture.  Nevertheless, the proof provided here 
is conceptually quite different from their argument (see Ferry-Weinberger [14] and
Hu [16] for related results on the Novikov conjecture).

Now let us further restrict to groups which are torsion-free $\delta$-hyperbolic.  For
such a group $\Gamma$, Theorem 1.1 above ensures that the group satisfies condition
$(*_\Delta)$.  In fact, $\delta$-hyperbolicity ensures that the $\Gamma$-action on
the pair $(\mathbb D^n, \Delta)$ has several additional properties.  In Section 4,
we will use these properties to show the following theorem:

\begin{Thm}
Let $\Gamma$ be a torsion-free $\delta$-hyperbolic group.  Then for each integer 
$n\geq 0$, the group homomorphism:
$$\bigoplus _{S\in \mathcal M} \pi_n(\phi_S):\bigoplus _{S\in \mathcal M} \pi_n
(\mathcal P (BS))\longrightarrow \pi_n(\mathcal P (B\Gamma))$$
is monic.
\end{Thm}

In the theorem above, $\mathcal M$ is a maximal collection of maximal infinite
cyclic subgroups of $\Gamma$, with no two elements in $\mathcal M$ being conjugate,
$\mathcal P (\cdot )$ is the stable topological pseudo-isotopy functor, and 
$\phi _S:\mathcal P (BS)\rightarrow \mathcal P (B\Gamma)$ is the functorially 
defined continuous map induced by $S\leq \Gamma$ (see Hatcher [15]).  We refer the
reader to section 4 for a more complete discussion of this result.

\section{EZ-structure implies condition $(*_\Delta)$}

Let us fix a discrete group $\Gamma$ with an EZ-structure $(\bar X,Z)$.
In this section we will provide a proof of Theorem 1.1.  In order to do this, we will
use the EZ-structure $(\bar X,Z)$ to build a new EZ-structure 
of the form $(\mathbb{D}^n, \Delta)$, where $\Delta$ is a closed subset of $\partial
\mathbb{D}^n =S^{n-1}$.  Let us start with a series
of lemmas that will allow us to make the structure of $\bar X-Z$ more suitable to 
our purposes.

\begin{Lem}[Reduction to a complex]
Let $\Gamma$ be a group with an EZ-structure $(\bar X,Z)$.  Then there is an EZ-structure
$(\tilde K\cup Z,Z)$, where $\tilde K$ is the universal cover of a finite
simplicial complex.
\end{Lem}

\begin{Prf}
We first observe that the hypotheses for an EZ-structure imply that the group $\Gamma$
is the fundamental group of an aspherical compact ANR, namely $(\bar X-Z)/\Gamma$.  
By a result of 
West [24], any compact ANR is homotopy equivalent to a compact polyhedra $K$.  In 
particular $K$ is a $K(\Gamma,1)$.  A result of Bestvina (Lemma 1.4 in [4]) now implies
that $(\tilde K \cup Z,Z)$ is an EZ-structure.
\end{Prf}

Our next step is to ``fatten'' $K$ so that it is a manifold with boundary.  In order to
do this, we embed (simplicially) $K$ into a high dimensional ($n\geq 5$) copy of 
$\mathbb{R}^n$, 
and let $W$ be a regular neighborhood of $K$.  Note that $W$ is a compact manifold with
boundary, and denote by $r:W\rightarrow K$ a retraction of $W$ onto $K$.  Let 
the retraction $\tilde r:
\tilde W\rightarrow \tilde K$ be the $\Gamma$-equivariant lift of $r$.

\begin{Lem}[Reduction to a manifold with boundary]
The pair $(\tilde W \cup Z,Z)$ is an EZ-structure for $\Gamma$.
\end{Lem}

\begin{Prf}
We follow the argument of Lemma 1.4 in Bestvina [4].  We start by taking the diagonal
embedding of $\tilde W$ in $(\tilde W\cup \infty)\times (\tilde K \cup Z)$.  The first 
factor is the
one point compactification of $\tilde W$, while the map into the second factor is given
by $\tilde r: \tilde W \rightarrow \tilde K \hookrightarrow \tilde K\cup Z$.  The topology
on $\tilde K \cup Z$ comes from taking the closure of the image of this diagonal 
embedding.  Lemma 1.3 in Bestvina [4] shows that this is a Z-structure.  Furthermore,
by construction, the action of $\Gamma$ on $\tilde W$ extends to an action of $\Gamma$
on $\tilde W\cup Z$.  Hence we have an EZ-structure.
\end{Prf}

An identical argument can be used to show the following:

\begin{Lem}[Doubling across the boundary]
Let $(N\cup Z, Z)$ be an EZ-structure on $\Gamma$, and assume that $N$ is a manifold 
(with or without boundary).  Denote by $\mathcal{N}$ the space $(N\times I)/\equiv$, 
where we collapse each $p\times I$, $p\in \partial N$, to a point (so if $N$ has no 
boundary, then $\mathcal{N}=N\times I$).  Then $(\mathcal{N}\cup Z,Z)$ is an EZ-structure
on $\Gamma$.
\end{Lem}

\begin{Prf}
We proceed as in the previous lemma, using the obvious $\Gamma$-equivariant map
$\rho: \mathcal{N}\rightarrow N\hookrightarrow N\cup Z$ in the place of $\tilde r$.
That is to say, we embed $\mathcal N$ into the space $(\mathcal N\cup \infty)\times
(N\cup Z)$ using the inclusion map on the first factor, and the map $\rho$ on the
second factor.  $\mathcal N \cup Z$ is then the closure of the image of $\mathcal N$ 
under this map, with the induced topology.  Once again, $Z$ lies as a Z-set, and the
mapping is $\Gamma$-equivariant by construction.
\end{Prf}

Note that the space $\mathcal{N}$ defined in Lemma 2.3 is also a 
manifold with boundary, and that the boundary $\partial \mathcal{N}$ of $\mathcal{N}$ 
is by construction
just the double of $N$ (the two copies being $N\times \{0\}$ and 
$N\times \{1\}$).

We now return to the situation we are interested in.  We have shown that we can reduce
to the case where the EZ-structure is of the form $(\tilde W\cup Z,Z)$, where 
$\tilde W$ is a manifold with boundary.  This allows us to apply the construction from 
the previous Lemma to obtain a new EZ-structure $(\mathcal{W}\cup Z, Z)$.  Our next result 
shows that $\mathcal{W}\cup Z$ is in fact a
topological manifold.  Because we will be refering to this result later in this section,
we prove it in a slightly more general form. 

\begin{Prop}
Let $(N\cup Z,Z)$ be an EZ-structure on $\Gamma$, and assume that $N$ is a manifold (with
or without boundary) of dimension $\geq 5$.  Let $(\mathcal N \cup Z, Z)$ be the 
EZ-structure defined in Lemma 2.3.  Then the space $\mathcal N \cup Z$ is a manifold 
with boundary.
\end{Prop}

\begin{Prf}
In order to show that the space $\mathcal N\cup Z$ is a compact manifold with 
we will use the celebrated characterization of high dimensional topological manifolds due 
to Edwards and Quinn (for a pleasant general survey, we refer to Mio [20]).  Recall that
this characterization provides a list of five necessary and sufficient conditions for a 
locally compact
high dimensional topological space to be a closed topological manifold.  The 
corresponding characterization for manifolds with boundary requires an additional 
condition about the `boundary'.  We will verify each of these six conditions as a 
separate claim.

\begin{Claim}[Finite dimensional]
The space $\mathcal N\cup Z$ is finite dimensional.
\end{Claim}

\begin{Prf}
Note that, by definition, $\mathcal N \cup Z$ is obtained by taking the closure of 
embedding of $\mathcal N$ into the space $(\mathcal N\cup \infty)\times (N\cup Z)$.
Both of the spaces $\mathcal N\cup \infty$ and $N\cup Z$ are finite dimensional, hence
so is their product $(\mathcal N \cup \infty)\times (N\cup Z)$.  Finally, $\mathcal N\cup
Z$ is a subset of a finite dimensional space, hence must also be finite dimensional.
\end{Prf}

\begin{Claim}[Locally contractible]
The space $\mathcal{N}\cup Z$ is locally contractible.
\end{Claim}

\begin{Prf}
This follows from the fact that the pair $(\mathcal N\cup Z,Z)$ is a Z-structure.  Indeed,
the first condition for a Z-structure forces $\mathcal N\cup Z$ to be an ER, and ER's 
are locally contractible.
\end{Prf}

\begin{Claim}[Homology manifold]
The space $\mathcal N\cup Z$ is a homology manifold with boundary.
\end{Claim}

\begin{Prf}
Let $n$ be the dimension of the manifold $N$.  We need to verify that the local homology 
of every point is either that of an $n$-dimensional sphere
(for ``interior'' points) or that of a point (for ``boundary'' points).  In order to do
this, we first observe that the local homology is easy to compute for points in 
$\mathcal N$.  Indeed, $\mathcal N$ is actually a manifold with boundary, hence the
local homology has the correct values.

Now let us focus on a point $p$ that lies on $Z \subset \mathcal N\cup Z$.  We claim 
that the (reduced) local homology at $p$ is trivial.  So we need to show that
$\bar H_*((\mathcal N\cup Z), (\mathcal N\cup Z)-p)=0$.  But this is also an immediate 
consequence of the fact that
$Z$ is a Z-set in $\mathcal N\cup Z$.  Indeed, an equivalent formulation of the Z-set
property states that there is a homotopy $J:(\mathcal N\cup Z)\times I\rightarrow 
\mathcal N\cup Z$
which satisfies the conditions:
\begin{itemize}
\item $J$ maps $\mathcal N\times I$ into $\mathcal N$.
\item $J_0:(\mathcal N\cup Z)\times \{0\}\rightarrow \mathcal N\cup Z$ is the identity map.
\item $J_t:(\mathcal N\cup Z)\times \{t\}\rightarrow \mathcal N\cup Z$ maps into 
$\mathcal N$ for all $t>0$.
\end{itemize}

In particular, the homotopy $J$ gives a family of homotopic maps which respect the
pair $((\mathcal N\cup Z),(\mathcal N\cup Z)-p)$, hence they all induce the same maps on 
the level of the homology groups $\bar H_*((\mathcal N\cup Z),(\mathcal N\cup Z)-p)$.  But 
the map induced by 
$J_0$ is the identity map, while the map induced by $J_1$ is the trivial map (since 
$J_1(\mathcal N\cup Z)\subset \mathcal N \subset (\mathcal N\cup Z)-p$).  Hence we have 
that the identity map coincides with the zero map, which immediately implies that
$\bar H_*((\mathcal N\cup Z),(\mathcal N\cup Z)-p)$ is trivial.
We conclude that $\mathcal N\cup Z$ is indeed a homology manifold with boundary.
\end{Prf}

Let us now recall the definition of the disjoint disk property.  A topological space 
$X$ has the disjoint disk property provided that any pair of maps from $\mathbb{D}^2$ 
into a space $X$ can be approximated, to an arbitrary degree of precision, by maps
whose images are disjoint.  

\begin{Claim}[Disjoint disk property]
The space $\mathcal N\cup Z$ has the disjoint disk property.
\end{Claim}

\begin{Prf}
Note that, since $\mathcal N\cup Z$ is an ER, it is metrizable; we will use this metric
to measure the closeness of maps.
Let $f,g$ be arbitrary maps from $\mathbb{D}^2$ into $\mathcal N\cup Z$, and 
let $\epsilon>0$ an arbitrary real number.  We need to exhibit a pair of maps which are
$\epsilon$ close to the maps we started with, and have disjoint image.

Observe that, since $Z$ is a Z-set in the space $\mathcal N \cup Z$, there is a map 
$H:\mathcal N\cup Z\rightarrow \mathcal N$ with the property that $H$ is an 
$(\epsilon /2)$-approximation of the identity map on $\mathcal N\cup Z$.
Consider the compositions $f^\prime :=H\circ f$ and $g^\prime :=H\circ g$, and observe 
that the maps $f^\prime$ and $g^\prime$ are $(\epsilon /2)$-approximations of $f$ and
$g$ respectively.  Furthermore, $f^\prime$ and $g^\prime$ map $\mathbb{D}^2$ into the
subset $\mathcal N$, which we know is a manifold 
of dimension $\geq 6$.  

But high dimensional manifolds automatically have the disjoint
disk property, so we can find $(\epsilon/2)$-approximations $f^{\prime \prime},
g^{\prime \prime}$ to the maps $f^\prime$, $g^\prime$ whose images are disjoint.
It is immediate from the triangle inequality that the $f^{\prime \prime},
g^{\prime \prime}$ satisfy our desired properties.  Hence the space 
$\mathcal N \cup Z$ has the disjoint disk property.
\end{Prf}

\begin{Claim}[Manifold point]
The space $\mathcal N\cup Z$ has a manifold point.
\end{Claim}

\begin{Prf}
By a manifold point, we mean a point with a neighborhood homeomorphic to some $\mathbb R^n$.
This is clear, since $\mathcal N$ is actually a topological manifold.
\end{Prf}

We now remind the reader of the 
characterization of high dimensional topological manifolds due to Edwards-Quinn 
([7],[22],[23]):

\begin{Thm}[Characterization of topological manifolds.] 
Let $X$ be a locally compact topological space, $n\geq 5$ an integer. 
Assume that $X$ satisfies the following properties:
\begin{itemize}
\item $X$ has the local homology of an $n$-dimensional manifold.
\item $X$ is locally contractible.
\item $X$ has finite (covering) dimension equal to $n$.
\item $X$ satisfies the disjoint disk property.
\end{itemize}
Then there is a locally defined invariant $I(X)\in 8\mathbb{Z}+1$ with the property that
$X$ is a topological manifold if and only if $I(X)=1$.
\end{Thm}

The corresponding theorem for a manifold with boundary requires an additional modification
of the first two conditions.  Namely, one needs to replace them with the following:
\begin{itemize}
\item every point $p\in X$ has either the local homology of an $(n-1)$-sphere, or that 
of a point.
\item the subset of points having the local homology of a point, denoted by 
$\partial_h(X)$ (the ``homological'' boundary), is a topological manifold of dimension $n-1$.
\end{itemize}
Under these two conditions, the Edwards-Quinn result implies that the space $X$ is a 
topological manifold with boundary (and the set $\partial_h(X)$ is the boundary of the
manifold $X$) if and only if the locally defined invariant $I(X)=1$ (see
Theorem 3.4.2 in Quinn [21]).

As such, we have reduced our theorem to showing the following:

\begin{Claim}
The set $\partial _h(\mathcal N\cup Z)$ is a compact manifold of dimension one lower
than the dimension of $\mathcal N$.
\end{Claim}

\begin{Prf}
By the proof of claim 3, we know exactly what the set $\partial _h(\mathcal N\cup Z)$
is.  Namely, it consists of the set $\partial \mathcal N \cup Z$.  Note that the set
$\partial \mathcal N$ is just the double of $N$ across it's boundary.  In particular,
$\partial _h(\mathcal N\cup Z)$ is obtained by taking two copies of $N \cup Z$,
and identifying the two copies of $\partial N\cup Z$.  

We now claim that $\partial N\cup Z$ is a Z-set in the space $N\cup Z$.  In order to show
this we need to exhibit a map $f_\epsilon:N\cup Z\rightarrow N\cup Z$ that is 
$\epsilon$-close to the identity, and has $f_\epsilon (N\cup Z)\subset N-\partial N$.  
Note that since $Z$ is a Z-set in $N\cup Z$, there is a map $g$ that is $(\epsilon/2)$-close
to the identity, and maps $N\cup Z$ into $N$.  Next, observe that since $N$ itself is a 
manifold with boundary, $\partial N$ is a Z-set in $N$, which implies the existence of 
a map $h:N\rightarrow N-\partial N$ which is $(\epsilon/2)$-close to the identity.  
Composing the two maps and using the triangle inequality gives us our desired claim.

So we see that $\partial _h(\mathcal N\cup Z)$ is obtained by doubling a Z-compactification
$N\cup Z$ of an open manifold $Int(N)$ along it's Z-boundary $\partial N\cup Z$.  By
a result of Ancel-Guilbault (Theorem 9 in [1]), this is automatically a manifold.
The dimension claim comes from the fact that $\partial _h(\mathcal N \cup Z)$ contains
$\partial \mathcal N$, hence must be a manifold of the same dimension as $\partial 
\mathcal N$, which is one less than the dimension of $\mathcal N$.
\end{Prf}

The Edwards-Quinn result now applies, completing our proof.
\end{Prf}

Let us summarize what we have so far: if $\Gamma$ has an EZ-structure,
we have shown that there is an EZ-structure $(\mathcal{W}\cup Z, Z)$ with the additional
property that $\mathcal W\cup Z$ is a topological manifold, and $Z$ is a closed subset
in the boundary of the topological manifold.  We now want to further improve the 
EZ-structure so that the space is in fact a topological disk.  In order to do this,
we iterate our procedure once more and define the space $\mathbb{W}=(\mathcal{W}\times
I)/\equiv$, where again the equivalence relation is given by collapsing $p\times I$,
$p\in \partial \mathcal W$ to a point.  By Lemma 2.3, the pair $(\mathbb{W}\cup Z,Z)$
is again an EZ-structure for $\Gamma$, and by Proposition 2.1, $\mathbb W\cup Z$ is 
a topological manifold with boundary.  We claim that $\mathbb W\cup Z$ is in fact a 
topological disk.

\begin{Prop}
The space $\mathbb W\cup Z$ is a disk.
\end{Prop}

\begin{Prf}
We begin by showing that the space $\partial (\mathbb W\cup Z)$ is simply connected.  
Notice that $\partial (\mathbb W\cup Z)$ is the double of the compact manifold with 
boundary $\mathcal W \cup Z$ along its boundary $\partial \mathcal W\cup Z$.  Furthermore
each of the spaces $\mathcal W\cup Z$ is contractible.  Siefert-Van Kampen now yields that
the double $\partial (\mathbb W\cup Z)$ must be simply connected.  Furthermore, observe
that the space $\mathbb W \cup Z$ is contractible.

Finally we note that any contractible manifold of dimension $\geq 6$ with simply 
connected boundary must
be homeomorphic to a disk.  This is a well known consequence of the h-cobordism theorem.  
A proof in the smooth
category can be found in Chapter 9, Proposition A, of Milnor's book [19].  
The same proof holds verbatim, replacing the use of Smale's smooth h-cobordism theorem
with the topological h-cobordism theorem of Kirby-Seibenmann's [18].
This concludes our proof of the proposition.
\end{Prf}

We have shown how given an arbitrary EZ-structure on a discrete group $\Gamma$, we can
construct an EZ-structure of the form $(\mathbb D^n, \Delta)$, where $\Delta$ is a
closed subset of $\partial \mathbb D^n=S^{n-1}$.  In particular, we see that any group
which has an EZ-structure automatically satisfies condition $(*_\Delta)$.

\section{Condition $(*_\Delta)$ implies the Novikov conjecture.}

We start this section by giving a reformulation of condition $(*_\Delta)$ which is closer to 
the formulation given by Farrell-Hsiang:

\begin{Def}
We say that a group $\Gamma$ satisfies condition $(*_\Delta)$ if for some integer $n$ there is an
action of $\Gamma$ on $(\mathbb{D}^n, \Delta)$, $\Delta$ a closed subset of $S^{n-1}=
\partial \mathbb{D}^n$ with the following two properties:
\begin{itemize}
\item $\Gamma$ acts properly discontinuously and cocompactly on $\mathbb{D}^n-\Delta$,
\item for each compact subset $K$ of $\mathbb{D}^n-\Delta$, and each $\epsilon>0$, there
exists a $\delta=\delta(K,\epsilon)>0$ such that for each $\gamma \in \Gamma$, if 
$d(\gamma K,\Delta)<\delta$, then $diam(\gamma K)<\epsilon$. 
\end{itemize}
\end{Def}

Observe that condition $(*_\Delta)$ generalizes condition (*) formulated
in Farrell-Hsiang [8] (the reader is also referred to [9] and the survey papers [12],[13]).
The only difference between the two conditions is that condition (*) also required
the set $\Delta$ to be $\partial \mathbb{D}^n=S^{n-1}$, and $\Gamma$ to be torsion-free.  
Furthermore, for torsion-free groups, it is easy to see that condition 
$(*_\Delta)$ corresponds exactly to the existence of an EZ-structure of the form $(\mathbb D^n,
\Delta)$, where $\Delta$ is a closed subset of $S^{n-1}$.

Note that, by the Theorem proved in the previous section, any group which has an EZ-structure
automatically satisfies condition $(*_\Delta)$.  In particular, the following two families
of groups satisfy condition $(*_\Delta)$:
\begin{itemize}
\item torsion-free $\delta$-hyperbolic groups.
\item torsion-free $CAT(0)$-groups.
\end{itemize}

Before starting the proof of Theorem 1.2, we first state the following useful Lemma:

\begin{Lem}
Let $(\mathbb{D}^m,\Delta)$ be a $\Gamma$-space satisfying the properties given in condition
$(*_\Delta)$.  Then there is a second $\Gamma$-space $(\mathbb{D}^{m+1}, \Delta)$ 
also satisfying $(*_\Delta)$, and a continuous
$\Gamma$-equivariant surjection $\mathbb{D}^m\times I\rightarrow \mathbb{D}^{m+1}$ mapping
$\Delta \times I$ to $\Delta$ and mapping $(\mathbb{D}^m-\Delta)\times I$ homeomorphically
to $\mathbb{D}^{m+1}-\Delta$.
\end{Lem}

\begin{Prf}
Let $\bar X =(\mathbb D^m \times I)/\equiv$, where the equivalence relation collapses each
line segment $x\times I$, $x\in \Delta$, to a point.  Let $\phi: \mathbb D^m \times I 
\rightarrow \bar X$ be the quotient map, and give $\bar X$ the $\Gamma$-space structure 
such that $\phi$ is $\Gamma$-equivariant.  Clearly, $\phi|_{(\mathbb D^m-\Delta)\times I}$
is a homeomorphism onto $\bar X -\Delta$.  

Projection onto the first factor of $\mathbb D^m\times I$ induces a $\Gamma$-equivariant
map $\Psi: \bar X-\Delta \rightarrow \mathbb D^m$.  The topology on $\bar X=(\bar X-
\Delta)\cup \Delta$ induced, using $\Psi$, by the construction in Lemma 2.2 coincides
with the one described above, as both topologies are compact and Hausdorff.  Hence
$(\bar X, \Delta)$ is an EZ-structure on $\Gamma$.

It remains to show that $\bar X$ is homeomorphic to $\mathbb{D}^{m+1}$.  For this we 
introduce a second decomposition space $Y=\mathbb D^m \times [0,2]/\sim$, where $\sim$
collapses each line segment $x\times [0,1]$, $x\in \Delta$, to a point.  Since $Y$ and
$\bar X$ are clearly homeomorphic, it suffices to construct a homeomorphism from 
$Y$ to $\mathbb D^m\times [0,2]$.  To do this, let $\phi :\mathbb D^m\rightarrow [0,1]$
be a continuous function such that $\phi^{-1}(0)=\Delta$.  Define 
$f:\mathbb D^m\times [0,2]\rightarrow \mathbb D^m\times [0,2]$ to be 
$f(x,t)=(x, t\phi (x))$ if $0\leq t\leq 1$, and
$f(x,t)=(x, (2-\phi (x))t + 2\phi(x)-2)$ if $1\leq t\leq 2$.  Observe that $f$ is a 
surjection.  

Since the point inverses of $f$ give the decomposition $\sim$ of $\mathbb D^m \times
[0,2]$, $f$ induces the desired homeomorphism.
\end{Prf}

The condition (*) was introduced by Farrell-Hsiang in order to provide an abstract setting
in which Novikov's Conjecture could be verified.  But the proof given in their paper 
carries over almost verbatim to the more general setting of condition $(*_\Delta)$.  Namely the 
following is true:

\begin{Thm}
Let $(\mathbb{D}^m,\Delta)$ be a $\Gamma$-space with the properties given in condition $(*_\Delta)$.
Suppose that $\Gamma$ is torsion-free, and let $M^m$ denote the orbit space $(\mathbb{D}^m
-\Delta)/\Gamma$.  Observe that $M^m$ is an aspherical compact manifold with boundary.  
Then the map in the (simple) surgery exact sequence:
$$\mathcal{S}^s(M^m\times \mathbb{D}^n,\partial)\longrightarrow [M^m\times \mathbb{D}^n,
\partial ; G/Top]$$
is identically zero when $n\geq 1$ and $n+m\geq 6$.
\end{Thm}

\begin{Prf}
For the reader's convenience, we recall the argument of [8] for the special case where
$\Gamma$ satisfies condition (*), as exposited in the Trieste notes [13], emphasizing the
modifications needed for the more general setting of condition $(*_\Delta)$.  So as not to obscure
the argument, we assume that $n=1$ and $M^m$ is triangulable.  Notice that the Lemma 3.1
formally reduces the general case $n\geq 1$ to the special case $n=1$.

Let $(\mathbb{D}^{m+1},\Delta)$ be the $\Gamma$-space determined by applying Lemma 3.1 to
the $\Gamma$-space $(\mathbb{D}^m,\Delta)$, and notice that $M^m\times \mathbb{D}^1=
(\mathbb{D}^{m+1}-\Delta)/\Gamma$.  Define the space:
$$\mathcal{E}^{2m+1}=(\mathbb{D}^{m+1}-\Delta)\times _\Gamma (\mathbb{D}^m-S^{m-1})$$
and let $p:\mathcal{E}^{2m+1}\rightarrow M^m\times \mathbb{D}^1$ be the bundle projection
induced by the projection to the {\it first} factor (the fiber of this projection is
$\mathbb{D}^m-S^{m-1}$).  Then the following diagram commutes:
$$
\begin{matrix}
\mathcal{S}^s(M^m\times \mathbb{D}^1,\partial) & \longrightarrow & [M^m\times \mathbb{D}^1,
\partial ; G/Top] \\
& & \\
\alpha \downarrow & & \downarrow p^*\\
& & \\
\mathcal{S} (\mathcal{E}, \partial) & \longrightarrow & [\mathcal{E},\partial; G/Top] \\
\end{matrix}
$$
where $\alpha$ is the obviously defined transfer map (see [13], pgs. 246-247).  Since 
$p$ is a homotopy equivalence, $p^*$ is an isomorphism.  Hence to prove the theorem, it 
is sufficient to verify the following:

\vskip 5pt

\noindent {\bf Assertion:} The map $\alpha$ is identically zero.

\vskip 5pt

To verify this assertion, note first that an arbitrary element in $\mathcal{S}^s
(M^m\times \mathbb{D}^1,\partial)$ can be represented by a pair $(f,h)$, where $f:M^m
\rightarrow M^m$ is a self-homeomorphism with $f|_{\partial M^m}=Id_{\partial M^m}$,
and $h:M^m\times \mathbb{D}^1 \rightarrow M^m \times \mathbb{D}^1$ is a homotopy of 
$f$ to $Id_{M^m}$ relative $\partial M^m$.  Define:
$$E^{2m}=(\mathbb{D}^{m}-\Delta)\times _\Gamma (\mathbb{D}^m-S^{m-1})$$
and notice that by Lemma 3.1, we have that $\mathcal{E}^{2m+1}=E^{2m}\times I$.

Observe that, given such a pair $(f,h)$, there is a well defined lift $\tilde f:
\mathbb{D}^m-\Delta \rightarrow \mathbb{D}^m-\Delta$, and that $\tilde f|_{S^{m-1}-\Delta}
=Id_{S^{m-1}-\Delta}$.  Now let $\tilde h$ be the unique lift of $h$ to 
$(\mathbb{D}^m-\Delta)\times I =\mathbb{D}^{m+1}-\Delta$ with the property that $\tilde h$
is a proper homotopy equivalence (relative $S^{m-1}-\Delta$) between $Id_{\mathbb{D}^m-
\Delta}$ and the self-homeomorphism $\tilde f$.

Then $k:=\tilde h\times Id_{\mathbb{D}^m-S^{m-1}}$ determines a proper homotopy (relative 
$\partial E$):
$$k:\mathcal{E}=E\times I \longrightarrow E\times I$$
between $Id_E$ and a self-homeomorphism $g:E\rightarrow E$ (which is also determined by
$\tilde f \times Id_{\mathbb{D}^m-S^{m-1}}$).  Note that $\mathcal{S} (\mathcal{E},\partial)
=\mathcal{S} (E\times I, \partial)$, since $\mathcal{E} =E\times I$.  Hence the pair
$(g,k)$ represents the image of the pair $(f,h)$ under the transfer map, i.e. $(g,k)=\alpha
(f,h)$.  The assertion then claims that the pair $(g,k)$ obtained in this manner is always 
zero in $\mathcal{S} (\mathcal{E} ,\partial)$.  In particular, the assertion would follow
from the following:

\begin{Prop}
$g$ is pseudo-isotopic to $Id_E$ (relative $\partial E$), via a pseudo-isotopy which is 
properly homotopic to $k$ (relative $\partial$).
\end{Prop}

We will now use the condition $(*_\Delta)$ to construct the pseudo-isotopy posited in this 
proposition.  Start by defining a new space $\bar E:= \mathbb{D}^m\times _\Gamma 
(\mathbb{D}^m-S^{m-1})$.  Note that the projection onto the {\it second} factor determines
a fiber bundle projection $q:\bar E\rightarrow Int(M^m)$ with fiber $\mathbb{D}^m$ (recall
that $Int(M^m)=(\mathbb{D}^m-S^{m-1})/\Gamma$).  Hence $\bar E$ is a manifold containing
$E$ as an open dense subset, and $\partial E\subset \partial \bar E$.  

Next observe that the second property of condition $(*_\Delta)$ implies that $\tilde f$ extends
to a $\Gamma$-equivariant homeomorphism $\bar f:\mathbb{D}^m\rightarrow \mathbb{D}^m$
by setting $\bar f|_{S^{m-1}} =Id _{S^{m-1}}$.  Consequently, $\bar f \times Id_{\mathbb{D}
^m-\Delta}$ determines a self-homeomorphism $\bar g:E\rightarrow E$ which extends 
$g:E\rightarrow E$ and satisfies $\bar g|_{\partial \bar E}=Id_{\partial \bar E}$.  We
now proceed to construct a pseudo-isotopy $\phi:\bar E\times I\rightarrow \bar E\times I$
satisfying:
\begin{itemize}
\item $\phi |_{\bar E\times \{0\}}=\bar g$
\item $\phi |_{\bar E\times \{1\}}=Id_{\bar E\times \{1\}}$
\item $\phi |_{(\partial \bar E)\times I}=Id_{(\partial \bar E)\times I}$
\end{itemize}
Once this is done, then the restriction of $\phi$ to the subset $E\times I \subset 
\bar E\times I$ will be the pseudo-isotopy posited in the proposition.

Observe that the three properties stated above define $\phi$ on the entire set $\partial
(\bar E\times I)$.  We need to extend $\phi$ over $Int(\bar E\times I)$.  In order to do 
this, consider the fiber bundle $r:\bar E\times I \rightarrow Int(M)$ with fiber 
$\mathbb{D}^m\times I$, where $r$ is the composition of the projection onto the first 
factor of $\bar E\times I$ followed by the map $q:\bar E\rightarrow Int(M)$.  Observe
that if $\sigma$ is an $n$-simplex in a triangulation of $Int(M)$, then $r^{-1}(\sigma)$
can be identified with $\mathbb{D}^{n+m+1}$.

The construction of $\phi$ proceeds by induction over the skeleta of $Int(M)$ via a 
standard obstruction theory argument.  And the obstructions encountered in extending 
$\phi$ from the $(n-1)$-skeleton to the $n$-skeleton are precisely those of extending a
self-homeomorphism of $S^{n+m}$ to a self-homeomorphism of $\mathbb{D}^{n+m+1}$.  But
these obstructions all vanish, because of the Alexander Trick.  Recall that this Trick
asserts that any self-homeomorphism $\eta$ of $S^n$ extends to a self-homeomorphism
$\bar \eta$ of $\mathbb{D}^{n+1}$.  In fact, $\bar \eta (tx)=t\eta (x)$ where $x\in S^n$
and $t\in I$ is an explicit extension.

Now the restriction $\psi := \phi|_{E\times I}$ is the pseudo-isotopy from $g$ to $Id_E$ 
asserted in the proposition.  A similar argument, which we omit, shows that $\psi$ is
properly homotopic to $k$ relative $\partial$.  This concludes the proof.
\end{Prf}

\section{Bounding $\pi_n(\mathcal P (B\Gamma))$ for $\delta$-hyperbolic groups.}

In this section, we give an application of our main result to obtaining a lower
bound for the homotopy groups $\pi_n(\mathcal P (B\Gamma))$ which holds for all
torsion-free
$\delta$-hyperbolic groups $\Gamma$.  Here $\mathcal P(\cdot )$ is the stable
topological pseudo-isotopy functor (see Hatcher [15]).  For this we need to first
recall some basic facts about $\delta$-hyperbolic groups.  Let $\Gamma$ be a torsion
free $\delta$-hyperbolic group (we exclude the case $\Gamma=\mathbb Z$).  
Then the following are true:

\begin{Fact}
If $S$ is an infinite cyclic subgroup of $\Gamma$, then there is a maximal infinite
cyclic subgroup containing $S$.  Furthermore this maximal subgroup is unique.
\end{Fact}

\begin{Fact}
If $C$ is a maximal infinite cyclic subgroup of $\Gamma$, then its normalizer is 
$C$ itself.
\end{Fact}

\begin{Fact}
If $S_1$ and $S_2$ are a pair of maximal infinite cyclic subgroups of $\Gamma$, and
$\{S_i^\pm\}\subset \partial ^\infty \Gamma$ are the corresponding pair of points in the
boundary at infinity, then either $S_1=S_2$
or $\{S_1^\pm\}\cap \{S_2^\pm \}= \emptyset$.
\end{Fact}

\begin{Fact}
If $S$ is a maximal infinite cyclic subgroup of $\Gamma$, then $\gamma\cdot S^-\neq S^+$ 
for all $\gamma \in \Gamma$.
\end{Fact}


We briefly explain why each of these facts holds.  The existence part of {\bf Fact 1} 
follows from Proposition 3.16 in Bridson-Haefliger (pg. 465 in [5]), 
while uniqueness follows from {\bf
Fact 3}.  For a maximal infinite cyclic subgroup, the normalizer coincides with the
centralizer.  If the element is not in the group itself, this would yield a pair
of commuting elements, giving a $\mathbb Z^2$ in $\Gamma$, which is impossible, giving
us {\bf Fact 2}.  {\bf Fact 3} follows from the proof of Theorem 3.20 in Bridson-Haefliger
(pg. 467 in [5]).  {\bf Fact 4} is an easy consequence of {\bf Fact 3}.  

Now fix a set $\mathcal M$ where the elements of $\mathcal M$ are maximal infinite 
cyclic subgroups of $\Gamma$ with each conjugacy class represented exactly once.
For each $S\in \mathcal M$, let $\phi_S:\mathcal P (BS)\rightarrow \mathcal P 
(B\Gamma)$ be the functorially defined continuous map (see Hatcher [15]).  
Note that $BS=S^1$ for each $S\in \mathcal M$.   
Theorem 1.3 that we are going to prove in this section 
states that, for each integer $n\geq 0$, the group homomorphism:
$$\bigoplus _{S\in \mathcal M} \pi_n(\phi_S):\bigoplus _{S\in \mathcal M} \pi_n
(\mathcal P (BS))\longrightarrow \pi_n(\mathcal P (B\Gamma))$$
is an injection.

Note that $\pi_0(\mathcal P(S^1))\cong \mathbb Z_2 \oplus \mathbb Z_2\oplus \cdots$,
where there are countably infinite number of $\mathbb Z_2$'s (see Igusa [17]).
Furthermore, the Isomorphism Conjecture for $\mathcal P(B\Gamma)$ formulated by
Farrell-Jones [11] is equivalent to the assertion that the homeomorphisms in 
Theorem 1.3 are all isomorphisms together with the assertion that the Whitehead
groups $Wh(\Gamma \times
\mathbb Z^n)$ vanish for all $n$.

Let us now proceed to prove Theorem 1.3.  By Theorem 1.1, we know that we have a 
sequence of EZ-structures $(\mathbb D^m, \partial ^\infty \Gamma)$, defined for all
sufficiently large $m$, such that $\Gamma$ acts on $\mathbb D^m$ by orientation 
preserving homeomorphisms, and $(\mathbb D^{m+1},\partial ^\infty \Gamma)=
(\mathbb D^m, \partial ^\infty \Gamma) \times I$ (i.e. is $\mathbb D^m \times I/\equiv$
where each interval $x\times I$, with $x\in S^{m-1}$, is collapsed to a point).  
Furthermore, each $S\in \mathcal M$ determines a pair of distinct points $S^+, S^- 
\in \partial ^\infty \Gamma$.  We start our argument by showing:

\begin{Claim2} 
$(\mathbb D^m, \{S^\pm\})$ is an EZ-structure for $S$.
\end{Claim2}

\begin{Prf}[Claim 1]
To see this claim, we first note that a closed subset of a Z-set is still a
Z-set, hence the pair $(\mathbb D^m,\{S^\pm\})$ satisfies the first two conditions
for an EZ-structure.  We also have, by restriction, an action of the group $S$
on $\mathbb D^m$.  Observe that the condition on the translates of compact sets 
forming a null sequence is inherited from the corresponding property for the 
the $\Gamma$-action.  So we are left with showing that the S-action on 
$\mathbb D^m-\{S^\pm\}$ is fixed point free, properly discontinuous, and cocompact.

To see that the S-action on $\mathbb D^m-\{S^\pm\}$ is fixed point free, 
we note that the $\Gamma$-action on $\mathbb D^m-\partial^\infty \Gamma$ is 
fixed point free, 
hence if the $S$-action has a fixed point, it must lie in the set $\partial 
^\infty \Gamma -\{S^\pm\}$.  But recall that the action of a $\delta$-hyperbolic 
group on it's boundary at infinity is hyperbolic.  More precisely, for every element
$g\in \Gamma$ ($g\neq 1$), we have a pair of fixed points $\{g^\pm\}\subset \partial 
^\infty 
\Gamma$ with the property that, for any compact set $C$ in 
$\partial ^\infty \Gamma -\{g^\pm\}$, and any open
sets $g^+\subset U^+$, $g^-\subset U^-$, there exists a positive integer $N$ 
such that:

\begin{itemize}
\item $g^n\cdot C\subset U^+$ for every $n\geq N$.
\item $g^{-n}\cdot C\subset U^-$ for every $n\geq N$.
\end{itemize}

In the particular case we are interested in, we have that $\{g^\pm\}=\{S^\pm\}$ for
every element $g\in S$ ($g\neq 1$).  Now assume that $p\in \partial ^\infty \Gamma 
-\{S^\pm\}$
is fixed by some element $g \in S$.  Then since $\partial ^\infty \Gamma$ is Hausdorff,
we can find a pair of open neighborhoods $U^\pm$ around the points $S^\pm$ which do 
not contain the given point $p$.  By hyperbolicity of the action, we have that some
high enough power of $g$ must take $p$ into $U^+$.  Hence $g$ cannot fix the point
$p$.

To see proper discontinuity of the action, we again restrict to looking at points
in $\partial ^\infty \Gamma -\{S^\pm\}$.  Indeed, since the $\Gamma$-action on 
$\mathbb D^m-\partial ^\infty \Gamma$ is properly discontinuous, so is the $S$-action
on $\mathbb D^m-\partial ^\infty \Gamma$.  So if proper discontinuity fails, it must
do so at some point $p\in \partial ^\infty \Gamma-\{S^\pm\}$.  But note that, by
hyperbolicity of the action (and as $\partial ^\infty \Gamma$ is Hausdorff), we
can find a triple of pairwise disjoint open sets $U^0, U^+, U^- \subset \partial 
^\infty \Gamma$ with $p\in U^0$, $S^+\in U^+$, $S^-\in U^-$, and with the property
that 
$g^n\cdot p\in U^+, g^{-n}\cdot p\in U^-$ for $n$ large enough ($g$ here refers to
a generator for the subgroup $S$).  In particular, there are only finitely many 
points of the form $g^i\cdot p$ which can lie in the open set $U^0$.  Note that the
topology on $\partial ^\infty \Gamma$ with respect to which the $S$-action is hyperbolic
coincides with the one induced on $\partial ^\infty \Gamma$ when viewed as a subset
of $\mathbb D^m$.  This gives us proper discontinuity of the $S$-action.

Finally, to see cocompactness, we need to exhibit a compact set in 
$\mathbb D^m-\{S^\pm\}$ whose $S$-translates cover $\mathbb D^m-\{S^\pm\}$.  Recall
that the $\Gamma$-action on $\mathbb D^m-\partial ^\infty \Gamma$ is cocompact, and
fix a compact fundamental domain $K_\Gamma$ for the $\Gamma$-action.  Now consider 
the Cayley graph $Cay(\Gamma)$ of the group $\Gamma$ with respect to a finite 
symmetric generating set.  Define the
set $T_S\subset \Gamma$ as follows: for each $S$-orbit of the $S$-action on $\Gamma$,
pick out an element that {\it minimizes} the distance in $Cay(\Gamma)$ from the 
$S$-orbit to the identity element (note that this choice might not be canonical).  
$T_S$ will consist of the union of one such element from each of
the $S$-orbits in $\Gamma$.  Now define the set $K_S$ to be the union of $T_S\cdot
K_\Gamma$ with the compact subset $C_\epsilon\subset \partial ^\infty 
\Gamma -\{S^\pm\}$, where $C_\epsilon$ is defined to be the complement of the open
$\epsilon$-neighborhood of $\{S^\pm\}$ for a sufficiently small $\epsilon$.  
We claim that $K_S$ is a compact fundamental domain 
for the $S$-action on $\mathbb D^m-\{S^\pm\}$, if $\epsilon$ is small 
enough.

We start by arguing that the $S$-translates of $K_S$ do indeed cover 
$\mathbb D^m-\{S^\pm\}$.  This is easy to see, as the $S$-translates of $T_S$ yield
the entire group $\Gamma$, and hence the union of the $S$-translates of $T_S\cdot 
K_\Gamma$ will coincide with $\Gamma \cdot K_\Gamma =\mathbb D^m-\partial ^\infty
\Gamma$.  On the other hand, the $S$-translates of $C_\epsilon$ will cover $\partial
^\infty \Gamma -\{S^\pm\}$, if $\epsilon$ is small enough, because the action 
of $g$ is uniformly continuous on 
$\partial ^\infty \Gamma$, fixes $S^\pm$, and $g^{\pm n}\cdot x \rightarrow S^\pm$
as $n\rightarrow \infty$ for every $x\in \partial ^\infty \Gamma -\{S^\pm\}$. This
gives us that $S\cdot K_S=\mathbb D^m-\{S^\pm\}$.

So we are left with showing that $K_S$ is compact in $\mathbb D^m-\{S^\pm\}$, 
provided $\epsilon$ is small enough.  We 
start by showing that $S^\pm \notin \overline{T_S\cdot K_\Gamma}$ (the overline
refers to the closure in $\mathbb D^m$).  So let us assume, 
by way of contradiction, that $S^+
\in \overline{T_S\cdot K_\Gamma}$ (the argument for $S^-$ is completely analogous).  
Then we can find a 
sequence of points $\{x_i\}$ in $T_S\cdot K_\Gamma$ with the property that 
$\lim x_i=S^+$.  Since each of these points lies in a translate of $K_\Gamma$,
we can consider instead the sequence $\{h_i\} \subset T_S$ of elements 
in the group $\Gamma$ having the property that $x_i\in h_i\cdot K_\Gamma$.  Now
the fact that $\lim x_i=S^+$ in $\mathbb D^m$ is equivalent to the fact that 
$\lim h_i=S^+$ in the Cayley graph $Cay(\Gamma)$.  

Recall that, in $Cay(\Gamma)$, saying that $\lim h_i=S^+$ implies that the sequence 
$\{h_i\}$
is within a uniformly bounded distance $D$ from the sequence $\{g^n\}$ (where $g$ is
the generator for $S$, and $n$ ranges over non-negative integers).  On the other 
hand, the definition of the set $T_S$ now forces the entire sequence $\{h_i\}$ to
lie within distance $D$ of the identity.  Indeed, if an element $h\in T_S$
lies within $D$ of an element $g^i\in S$, then (since $\Gamma$ acts by isometries
on its Cayley graph) the element $g^{-i}h$ lies within $D$ of the identity element,
and is in the same $S$-orbit as the element $h$.  In particular, this forces $h$ to
be within $D$ of the identity (since by construction, $h$ minimizes the distance to 
the identity within its $S$-orbit).  So we have exhibited a bounded sequence in 
$Cay(\Gamma)$ which converges to a point in $\partial ^\infty \Gamma$, giving us 
our contradiction.  

So the closure of $T_S\cdot K_\Gamma$ does not contain $S^\pm$, hence the intersection
of the closure of $T_S\cdot K_\Gamma$ with $\partial ^\infty \Gamma$ lies outside of
a small $\epsilon$-neighborhood of $\{S^\pm\}$.  This immediately gives us that the
union of $T_S\cdot K_\Gamma$ with the corresponding $C_\epsilon$ is a compact set in
$\mathbb D^m-\{S^\pm\}$, and hence that the action of $S$ on $\mathbb D^m-\{S^\pm\}$ 
is cocompact.

We now have that the pair $(\mathbb D^m,\{S^\pm\})$ satisfies all the conditions
for an EZ-structure, concluding the proof of Claim 1.
\end{Prf}

We now continue the proof of Theorem 1.3.  Note that $(\mathbb D^{m+1}, \{S^\pm\}) 
=(\mathbb D^m, \{S^\pm\})\times I$.  Arguing as in the paper by Farrell-Jones 
(see pgs. 462-467 in [10]), it suffices 
to construct, for each sufficiently large integer $m$, a pair of continuous maps:
$$ 
\begin{matrix}
g_S :P (M^m_S) \longrightarrow P(M^m)\\
g^S :P (M^m) \longrightarrow P(M^m_S)\\
\end{matrix}
$$
where $M^m=(\mathbb D^m -\partial ^\infty \Gamma)/\Gamma$, $M^m_S=(\mathbb D^m-\{S^\pm
\})/S$, and $P(\cdot )$ denotes the (unstable) pseudo-isotopy space, and where the 
maps $g_S$ and $g^S$ satisfy the following:

\vskip 5pt

\noindent{\bf{Assertion:}} $g^S\circ g_S$ is homotopic to the identity, and 
$g^{S^\prime}\circ g_S$ is homotopic to a constant map whenever $S\neq S^\prime$.

\vskip 5pt

We first discuss the construction of the maps $g_S$, $g^S$, and will then discuss why
the pair of maps we constructed satisfy the
assertion.  Start by observing that both $M^m$ and $M^m_S$ are compact $m$-dimensional
manifolds with boundary (we will henceforth suppress the superscript indicating 
dimension unless it is explicitly relevant to the argument being presented).
Now let $p= p_S: Int(M_S)\rightarrow Int(M)$ be the covering space corresponding 
to the subgroup $S\subset \Gamma= \pi_1(Int(M))$.  Using the
s-cobordism theorem (and assuming $m\geq 6$), one easily constructs an isotopy 
$\phi _t =\phi^S_t:M_S\rightarrow M_S$ such that $\phi_0=Id_{M_S}$, and $p\circ \phi_1
:M_S\rightarrow M$ is an embedding.  To define $g_S$, let $f:M_S\times I\rightarrow
M_S\times I$ be a pseudo-isotopy (i.e. an element of $P(M_S)$).  Recall that $f$
is an automorphism (i.e. an onto homeomorphism) with the property that:
$$f|_{M_S\times \{0\} \cup (\partial M_S)\times I}= Id|_{M_S\times \{0\} \cup 
(\partial M_S)\times I}.$$

We can now define $f_*=g_S(f)\in P(M)$ by setting $f_*(x,t)$ to be:
\begin{itemize}
\item $(x,t)$ if $x\in M-Image(p\circ \phi_1)$
\item $p\circ \phi_1(f(\bar x, t))$ if $x=p\circ \phi_1(\bar x)$
\end{itemize}
where $x\in M$ and $t\in I$.  This gives us the map $g_S$.

On the other hand, to define $g^S(f)$, where $f\in P(M)$, let $\tilde f:
(\mathbb D^m-\partial ^\infty \Gamma)\times I\rightarrow (\mathbb D^m-\partial 
^\infty \Gamma)\times I$ be the lift of $f$ such that $\tilde f (x,t)=(x,t)$
if either $x\in S^{m-1}=\partial \mathbb D^m$ or if $t=0$.  Now $\tilde f$ induces
an automorphism $\bar f$ of $(\mathbb D^{m+1}, \partial ^\infty \Gamma)$, since
$(\mathbb D^{m+1}, \partial ^\infty \Gamma)=(\mathbb D^{m}, \partial ^\infty \Gamma)
\times I$.  Note that $\bar f$ is $\Gamma$-equivariant and that $\bar f|_{\partial
_- \mathbb D^{m+1}} = Id_{\partial_- \mathbb D^{m+1}}$, where ${\partial
_- \mathbb D^{m+1}}$ is the image of $\mathbb D^m\times \{0\} \cup S^{m-1}\times I$
under the quotient map $\mathbb D^m\times I\rightarrow \mathbb D^{m+1}$.
Since $\partial ^\infty \Gamma\subset {\partial_- \mathbb D^{m+1}}$, $\bar f$ 
induces an $S$-equivariant automorphism of $\mathbb D^{m+1}-\{S^\pm\}$ which then
descends to an automorphism $f_S$ of $(\mathbb D^{m+1}-\{S^\pm\})/S$.  After 
``appropriately identifying'' 
$$\mathcal M_S=(\mathbb D^{m+1}-\{S^\pm\})/S$$
with $M^m_S\times I$, $g^S(f)$ is defined by $g^S(f)=f_S$.  

To do this identification, first note that $\mathcal M_S$ is the quotient space of
$M^m_S\times I$ where each interval $x\times I$, $x\in \partial M^m_S$ is collapsed
to a point.  So $M^m_S\times \{0\}$ is canonically identified with a codimension
zero submanifold $\partial _-\mathcal M_S$ of $\partial \mathcal M_S$.  By equating
$\partial M^m_S\times I$ with a short collar of $\partial (\partial _-\mathcal M_S)$
in $\partial \mathcal M_S$, an identification of $M_S\times I$ to $\mathcal M_S$
can be constructed such that the composition:
$$P(M_S)\longrightarrow Aut(\mathcal M_S, \partial _-(\mathcal M_S))\longrightarrow
P(M_S)$$
is homotopic to the identity (here the two maps above are the naturally defined
continuous maps; in fact, the second map is the homeomorphism induced by the identification
while the first is determined by the fact that $\mathcal M_S$ is a quotient space of 
$M_S\times I$).  This is the ``appropriate identification'' mentioned above.

This gives us the two maps for which we claim the assertion holds.  Before continuing
our proof, we note that, when $m\geq 6$, the spaces $M^m_S$ are all homeomorphic to
$S^1\times \mathbb D^{m-1}$.  Indeed, this follows by the s-cobordism theorem, and the
fact that $S$ acts via orientation preserving homeomorphisms on $\mathbb D^m-\{S^\pm\}$;
thus the closed tubular neighborhood of any embedded circle $S^1$ in $Int(M^m_S)$,
which induces a homotopy equivalence, is homeomorphic to $S^1\times \mathbb D^{m-1}$.

Now the {\bf Assertion}, made above, can be verified in the same way that properties
(i) and (ii) in Lemma 2.1 of Farrell-Jones [10] were proven.  We merely point out 
that they follow directly from the following two claims which we proceed to formulate
and then to verify.  Let $T_S$ denote the image of $p_S\circ \phi^S_1$.  Note that
$T_S$ is a codimension zero submanifold of $Int(M^m_S)$ and that $T_S$ is homeomorphic
to $S^1\times \mathbb D^{m-1}$.  Recall that:
$$p_S:Int(M_S)\longrightarrow Int(M)$$
is the covering projection corresponding to $S\subset \Gamma$.  And that $\phi^S_1:
M_S\rightarrow Int(M_S)$ is an embedding isotopic to $Id_{M_S}$.  Recall that we  
assumed that $\Gamma$ is not cyclic.

Now let $\{C_i\}$ denote the connected components of $p_S^{-1}(T_S)$, and note that
$p_S^{-1}=\coprod _i C_i$.  Let $\bar C_i$ denote the closure of $C_i$ in $M_S$.  It
is an elementary observation that each $C_i$ is a codimension zero submanifold of
$Int(M_S)$ as well as an open subset of $p_S^{-1}(T_S)$.  Furthermore, observe that
$Image(\phi_1^S)$ is a codimension zero submanifold of $Int(M_S)$ which is homeomorphic
to $S^1\times \mathbb D^{m+1}$.

\begin{Claim2}
We can index the set $\{C_i\}$ so that $C_0=Image(\phi^1_S)$ and $\bar C_i$ is
homeomorphic to $\mathbb D^{m}$ when $i\neq 0$.
\end{Claim2}

Now let $S^\prime \in \mathcal M$ with $S^\prime \neq S$, and denote by $\{C_i^\prime\}$
the connected components of $p_{S^\prime}^{-1}(T_S)$ and by $\bar C_i^\prime$ the
closure of $C^\prime _i$ in $M_{S^\prime}$.  It is again elementary that each 
$C_i^\prime$ is a codimension zero submanifold of $Int(M_{S^\prime})$ as well as an 
open subset of $p_{S^\prime}^{-1}(T_S)$.

\begin{Claim2}
Each $\bar C_i^\prime$ is homeomorphic to $\mathbb D^m$.
\end{Claim2}

We now proceed with the proofs of the two claims.  The Facts $1_\delta$-$4_\delta$ 
used in the proofs below refer 
to the facts about $\delta$-hyperbolic groups discussed at the beginning of this 
section. 

\begin{Prf}[Claim 2]
One easily sees that each $p_i:C_i\rightarrow T_S$ is a covering projection where
$p_i=p_S|_{C_i}$.  Hence $Image(\phi_1^S)$ must be one of the components $C_i$ since
$p:Image(\phi_i^S)\rightarrow T_S$ is a homeomorphism.  Thus we may index the 
components starting with $C_0=Image(\phi_1^S)$.  Therefore it remains to show that 
$\bar C_i$ is homeomorphic to $\mathbb D^m$ when $i\neq 0$.  To do this, define
\begin{itemize}
\item $q:\mathbb D^m-\partial ^\infty \Gamma \longrightarrow M=(\mathbb D^m-\partial
^\infty \Gamma)/\Gamma$
\item $r=r_S:\mathbb D^m-\{S^\pm\} \longrightarrow M_S=(\mathbb D^m -\{S^\pm\})/S$
\end{itemize}
to be the universal covering maps whose groups of deck transformations are $\Gamma$
and $S$ respectively.  Then we have the following commutative triangle of covering
spaces:

$$\xymatrix{ Int(\mathbb D^m) \ar[rr]^r \ar[dr]_q & & Int(M_s) \ar[dl]^p
\\ & Int(M) &
}
$$

Note that $q^{-1}(T_S)=\coprod _i D_i$ where each $D_i$ is a connected component of
$q^{-1}(T_S)$.  And let $\bar D_i$ be the closure of $D_i$ in $\mathbb D^m$.  One 
easily sees the following ten points:
\begin{enumerate}
\item Each $D_i$ is open in $q^{-1}(T_S)$.
\item Each $D_i$ is a codimension zero submanifold of $Int(\mathbb D^m)$.\
\item $q_i:D_i\rightarrow T_S$ is a universal covering space (where $q_i=q|_{D_i}$)
whose group of deck transformations $S_i$ consists of all $\gamma \in \Gamma$ such
that $\gamma (D_i)=D_i$.  Consequently, $D_i$ is homeomorphic to $\mathbb D^{m-1}\times
\mathbb R$.
\item The components $D_i$ are permuted transitively by $\Gamma$.  Consequently, the
groups $S_i$ are all conjugate cyclic subgroups of $\Gamma$.
\item At least one of the groups $S_i$ is $S$.  Hence all the $S_i$ are maximal cyclic
subgroups of $\Gamma$.  And we can rearrange the indexing so that $S_0=S$.
\item If the cardinality $|S_i\cap S_j|>1$, then $i=j$.  This follows from points
(4) and (5) by using Fact $1_\delta$ and Fact $2_\delta$.
\item Let $\tilde \phi _t:\mathbb D^m-\{S^\pm\}\rightarrow \mathbb D^m-\{S^\pm\}$ be
the lift of the isotopy $\phi _t$ with respect to the covering projection $r$ such
that $\tilde \phi_0= Id$.  Then $D_0=Image(\tilde \phi_1)$, and consequently
$\bar D_0=D_0\cup \{S^\pm\}$, which forces $\bar D_0$ to be homeomorphic to 
$\mathbb D^m$.
\item Because of points (4) and (7), $\bar D_i=D_i\cup \{S_i^\pm\}$ and is homeomorphic
to $\mathbb D^m$.  Also because of point (6) and Fact $3_\delta$, $\bar D_i\subset
\mathbb D^m-\{S^\pm\}$ if $i\neq 0$, and consequently $\bar D_i$ is also the closure
of $D_i$ in $\mathbb D^m-\{S^\pm\}$.
\item If $\gamma (\bar D_i)\cap \bar D_i\neq \emptyset$, where $\gamma \in \Gamma$, 
then $\gamma \in S_i$.  This results from points (4), (6), (8), along with Facts 
$3_\delta$ and $4_\delta$.  Consequently, if $i\neq 0$, then $r|_{\bar D_i}: \bar D_i
\rightarrow
r(\bar D_i)=\overline {r(D_i)}$ is a homeomorphism since $S_i\cap S_0=1$, because of point
(6) (Here $\overline {r(D_i)}$ denotes the closure of $r(D_i)$ in $M_S$).
\item There is a surjection of indexing sets $i\mapsto \alpha (i)$, with $\alpha (0)=0$,
such that $r_i:D_i\rightarrow C_{\alpha (i)}$ is a covering space (here $r_i$ denotes
$r|_{D_i}$).  This follows from the above commutative triangle in which $p$, $q$,
and $r$ are open maps.  
\end{enumerate}
It now follows immediately from points (8), (9), and (10), that $\bar C_i$ is 
homeomorphic to $\mathbb D^m$ when $i\neq 0$; thus completing the proof of Claim 2.
\end{Prf}

\begin{Prf}[Claim 3]
This proof closely parallels the one just given for Claim 2.  Note that the above 
points (1)-(9) continue to hold.  And by replacing $S$ by $S^\prime$ in the above 
commutative triangle, the following analogue (10)$^\prime$ of point (10) is similarly 
verified using that $p_{S^\prime}$, $q$, and $r_{S^\prime}$ are open maps: there is a 
surjection $i\mapsto \beta (i)$ of indexing sets such that $r_i^\prime : D_i
\rightarrow C_{\beta (i)}^\prime$ is a covering space where $r_i^\prime=r_{S^\prime}|
_{D_i}$.  

Then Fact $3_\delta$ yields that:
$$\{S_i^\pm\}\subseteq (\mathbb D^m-\partial ^\infty S^\prime)=Domain(r_{S^\prime})$$
which together with point (8) shows that 
$$\bar D_i \subseteq Domain(r_{S^\prime}).$$
Therefore point (9) yields that:
$$r_{S^\prime}|_{\bar D_i}:\bar D_i\longrightarrow r_{S^\prime}(\bar D_i)=\overline {r_
{S^\prime}(D_i)}=\overline {C_{\beta (i)}^\prime}$$
is a homeomorphism.  But $\bar D_i$ is homeomorphic to $\mathbb D^m$ by point (8),
and $\beta$ is a surjection by point (10)$^\prime$.  This concludes the proof of 
Claim 3.
\end{Prf}

Finally, we point out that, from these two claims, it is easy to show the assertion.
Indeed, the pseudo-isotopies $g^S\circ g_S(f)$ and $g^{S^\prime}\circ g_S(f)$ are
supported over $\cup_i\bar C_i$ and $\cup _i\bar C_i^\prime$ respectively.  Because of
claims 2 and 3, the Alexander trick can be used to verify the {\bf Assertion}.  We 
refer the reader to section 2 of Farrell-Jones [10] for more details.

\section{Concluding remarks.}

We would like to conclude by asking the question: which finitely generated 
groups have an EZ-structure?
A version of this question was already posed by Bestvina [3], where he asks whether
every group $\Gamma$ with a finite $B \Gamma$ has a Z-structure.  
It is interesting to construct groups which are neither $\delta$-hyperbolic, nor
$CAT(0)$ groups, but do have an EZ-structure.  Bestvina gives some important examples
of such groups in [3].  Do torsion free
subgroups of finite index in $SL_n(\mathbb Z)$ have an EZ-structure?

It would also be of some interest to find applications of Theorem 1.1 to geometric group
theory.  Indeed, condition $(*_\Delta)$ for torsion free groups yields an action of the 
group on disks, which, aside from a ``bad limit set'' is properly discontinuous, fixed 
point free, and cocompact.  With the exception of cocompactness, this is reminiscent
of the action of a Kleinian group
on (the compactification) of hyperbolic $n$-space.  In some sense, Theorem 1.1 states that
every torsion-free $\delta$-hyperbolic group has an action that mimics that
of a Kleinian group.  One feels that this should have some strong geometric consequences.

In addition, one could consider the possibility of strengthening condition $(*_\Delta)$
by also requiring the action of the group $\Gamma$ on $\mathbb D^n$ to be {\it smooth}.
Work of Benoist-Foulon-Labourie [2] suggests that among $\delta$-hyperbolic groups, 
perhaps only uniform
lattices satisfy this extra property.  In any event it would be interesting to determine
which $\delta$-hyperbolic groups satisfy this smooth form of condition $(*_\Delta)$.

\section{Bibliography}

\vskip 10pt

\noindent [1] Ancel, F.D. \& Guilbault, C.R.  {\it $\mathcal Z$-compactifications of open 
manifolds}.  Topology  38  (1999),  no. 6, pp. 1265--1280.

\vskip 5pt

\noindent [2] Benoist, Y., Foulon, P. \& Labourie, F. {\it Flots d'Anosov à distributions 
stable et instable différentiables}. (French) [Anosov flows with stable and unstable 
differentiable distributions]. J. Amer. Math. Soc. 5 (1992), no. 1, pp. 33--74.

\vskip 5pt

\noindent [3] Bestvina, M. \&  Mess, G. {\it The boundary of negatively curved groups}.
 J. Amer. Math. Soc.  4  (1991),  no. 3, pp. 469--481.

\vskip 5pt

\noindent [4] Bestvina, M.  {\it Local homology properties of boundaries of group}.
 Michigan Math. J.  43  (1996),  no. 1, pp. 123--139.

\vskip 5pt

\noindent [5] Bridson, M. \& Haefliger, A. {\it Metric spaces of non-positive curvature}.
Grundlehren der Mathematischen Wissenschaften [Fundamental Principles of Mathematical 
Sciences], 319. Springer-Verlag, Berlin,  1999. xxii+643 pp. 

\vskip 5pt

\noindent [6] Carlsson, G. \& Pedersen, Erik K.  {\it Controlled algebra and the 
Novikov conjectures for $K$- and  $L$-theory}.  Topology  34  (1995),  no. 3, 
pp. 731--758.

\vskip 5pt

\noindent [7] Edwards, R.D.  {\it The topology of manifolds and cell-like maps}. in
 Proceedings of the International Congress of Mathematicians (Helsinki, 1978), 
 pp. 111--127, Acad. Sci. Fennica, Helsinki,  1980. 

\vskip 5pt

\noindent [8] Farrell, F.T. \& Hsiang, W.C.  {\it On Novikov's conjecture for 
nonpositively curved manifolds. I}.  Ann. of Math. (2)  113  (1981),  no. 1,
pp. 199--209.

\vskip 5pt

\noindent [9] Farrell, F.T. \& Hsiang, W.C.  {\it On Novikov's conjecture for cocompact 
discrete subgroups of a Lie  group}. in Algebraic topology, Aarhus 1982 (Aarhus, 1982), 
pp. 38--48, Lecture Notes in Math., 1051, Springer-Verlag, Berlin,  1984. 

\vskip 5pt

\noindent [10] Farrell, F.T. \& Jones, L.E.  {\it $K$-theory and dynamics. II}.
 Ann. of Math. (2)  126  (1987),  no. 3, pp. 451--493.

\vskip 5pt

\noindent [11] Farrell, F.T. \& Jones, L.E.  {\it Isomorphism conjectures in algebraic 
$K$-theory}.  J. Amer. Math. Soc.  6  (1993),  no. 2, pp. 249--297.

\vskip 5pt

\noindent [12] Farrell, F.T.  {\it Lectures on surgical methods in rigidity}.
Published for the Tata Institute of Fundamental Research, Bombay; by  Springer-Verlag, 
Berlin,  1996. iv+98 pp. 

\vskip 5pt

\noindent [13] Farrell, F.T.  {\it The Borel conjecture}.  in Topology of high dimensional
manifolds, Vol 1 (eds. F.T. Farrell, L. Goettsche and W. Lueck), ICTP Lecture Notes, Trieste,
2002.

\vskip 5pt

\noindent [14] Ferry, S.C. \& Weinberger, S.  {\it Curvature, tangentiality, and controlled 
topology}.  Invent. Math.  105  (1991),  no. 2, pp. 401--414.

\vskip 5pt

\noindent [15] Hatcher, A.E.  {\it Concordance spaces, higher simple-homotopy theory, and  
applications}. in Algebraic and geometric topology (Proc. Sympos. Pure Math., Stanford  
Univ., Stanford, Calif., 1976), Part 1, pp. 3--21, Proc. Sympos. Pure Math., XXXII, 
Amer. Math. Soc., Providence, R.I.,  1978. 

\vskip 5pt

\noindent [16] Hu, B.  {\it Retractions of closed manifolds with nonpositive curvature}. in
Geometric group theory (Columbus, OH, 1992), pp. 135--147, Ohio State Univ. Math. Res. 
Inst. Publ., 3, de Gruyter, Berlin,  1995. 

\vskip 5pt

\noindent [17] Igusa, K.  {\it What happens to Hatcher and Wagoner's formulas for 
$\pi \sb{0}C(M)$  when the first Postnikov invariant of $M$ is nontrivial?} in 
Algebraic $K$-theory, number theory, geometry and analysis (Bielefeld,  1982), 
pp. 104--172, Lecture Notes in Math., 1046, Springer, Berlin,  1984. 

\vskip 5pt

\noindent [18] Kirby, R.C. \& Siebenmann, L.C.  {\it Foundational essays on topological 
manifolds, smoothings, and triangulations}.
Ann. of Math. Stud., 88,
Princeton University Press, Princeton, NJ, 1977. vii+355 pp.

\vskip 5pt

\noindent [19] Milnor, J.  {\it Lectures on the $h$-cobordism theorem}.
Princeton University Press, Princeton, NJ,  1965. v+116 pp.

\vskip 5pt

\noindent [20] Mio, W.  {\it Homology manifolds}. in
 Surveys on surgery theory, Vol. 1, 
 pp. 323--343, Ann. of Math. Stud., 145, Princeton Univ. Press, Princeton, NJ,  2000. 

\vskip 5pt

\noindent [21] Quinn, F.  {\it Ends of maps, I}.  Annals of Math.  110  (1979), pp. 
275--331.

\vskip 5pt

\noindent [22] Quinn, F.  {\it Resolutions of homology manifolds, and the topological 
characterization of manifolds}. Invent. Math.  72  (1983),  no. 2, pp. 267--284.

\vskip 5pt

\noindent [23] Quinn, F.  {\it An obstruction to the resolution of homology manifolds}.
 Michigan Math. J.  34  (1987),  no. 2, pp. 285--291.

\vskip 5pt

\noindent [24] West, J.E.  {\it Mapping Hilbert cube manifolds to ANR's: a solution of a 
conjecture of  Borsuk}.  Ann. Math. (2)  106  (1977),  no. 1, pp. 1--18.

\end{document}